\input amstex

\documentstyle{amsppt}

\topmatter

\title
On Double Cosets in Free Groups
\endtitle

\author
Rita Gitik and Eliyahu Rips
\endauthor

\address 
Institute of Mathematics, Hebrew University, Jerusalem, 91904, Israel
\endaddress
\curraddr 
A and H consultants, Ann Arbor, MI 48104
\endcurraddr
\email 
ritagtk \@ umich.edu
\endemail

\address 
Institute of Mathematics, Hebrew University, Jerusalem, 91904, Israel
\endaddress

\abstract
It is shown that for any finitely generated subgroups $H$ and $K$ 
of a free group $F$, and for any $g \in F$ the double coset $H gK$ is 
closed in the profinite topology of $F$.
\endabstract
\endtopmatter
\document

A well-known theorem of M. Hall \cite{2} states (in  different language) 
that any finitely generated subgroup of a free group is closed 
in the profinite topology. We show that a slight modification of its proof
\cite{6} yields a stronger result:

\proclaim{Theorem} For any finitely generated subgroups $H$ and $K$ of a free
group $F$, and for any $g \in F$ the double coset $H gK$ is closed in 
the profinite topology of $F$.
\endproclaim
A free group $F=<X>$ can be viewed as the fundamental group of a wedge 
$W$ of $|X|$ oriented circles labeled by elements of $X$. Subgroups of $F$ 
correspond
bijectively to based covering spaces of $W$, and any covering space
of $W$ is a graph which inherits orientation and labeling 
of its edges from $W$.

Let $X_0 \subset X$ and let $\Gamma$ be a subgraph of a covering of $W$. 
An $X_0$-component of $\Gamma$ is a maximal connected subgraph of $\Gamma$ 
with all its edges labeled by elements of $X_0$.

\proclaim{Lemma} Let $\Gamma$ be a subgraph of a covering of $W$ such that $\Gamma$ has 
finitely  many vertices. There exists an embedding of $\Gamma$ in a covering
$\Gamma'$ of $W$ such that $\Gamma$ and $\Gamma'$ have the same vertices, and
for any $X_0 \subset X$ distinct $X_0$-components of $\Gamma$ remain distinct
in $\Gamma'$.
\endproclaim

\demo{Proof}  We give an algorithm for constructing $\Gamma'$  by adding edges to
$\Gamma$ in a unique way. For any vertex $v$ of $\Gamma$ and for any $x \in X$
the number of edges labeled with $x$ having an endpoint at $v$ is either
$0,1$ or $2$. If the number is $0$, we add  an edge labeled with $x$
with both endpoints at $v$. 
If the number is $1$ or $2$, let $p_x$
be the maximal path consisting of edges labeled only with $x$ and with an
endpoint at $v$. If $p_x$ has  both endpoints at $v$ we do nothing. 
Otherwise we add to $\Gamma$ an edge labeled 
with $x$  connecting the endpoints of $p_x$. It is clear that the 
projection from $\Gamma$ to $W$ extends uniquely to a covering map from
$\Gamma'$ to $W$.
\enddemo

\remark{Remark}
 Let $H$ be a finitely  generated  subgroup of $F$, and let $f \in F 
\setminus H$. Let $\Gamma$
be the minimal connected subgraph of the covering $C$ of $W$ corresponding 
to $H$ which contains the core of $C$ (cf. \cite{6}) and the path $p$  
beginning at the basepoint $v_0$ of $C$ whose projection in $W$ represents $f$.
Embed $\Gamma$ in a covering $\Gamma'$ as in the lemma. Then as $\Gamma$  has finitely
many vertices, so does $\Gamma'$, therefore the subgroup $M$ of $F$ corresponding 
to $\Gamma'$ has finite index in $F$. As  $p$ is not a closed path in $\Gamma$,
it remains not closed in $\Gamma'$, hence $f \notin M$. As $\Gamma$ is a subgraph
of $\Gamma',\;\; H$ is a subgroup of  $M$, proving M. Hall's theorem (cf. 
\cite{6}).
\endremark

\demo{Proof of the theorem} As $HgK= H(gKg^{-1})g$, it is enough to consider the 
case $g=1$.

By an observation due to P.Kropholler, we can replace $F$ by a subgroup of 
finite  index (cf. \cite{3}), so we can assume  that $K$ is a free factor of
$F, \;\; F=K*L$.  Let $X_1$ and $X_2$ be sets of free generators of $K$ and $L$
respectively, then $X= X_1 \cup X_2$ is a set of free generators of $F$.
Let $f\in F, f \notin  HK$. Our goal is to construct  a subgroup $M$
of  finite index in $F$ such that $Mf \cap HK = \emptyset $.

Let $\Gamma, \Gamma'$ and $M$ be as in the remark. 
As $f \notin HK$, the $X_1$-component of $v_0$ in $\Gamma$ does not contain the
endpoint of $p$,
therefore the lemma implies that $\Gamma'$ has the same property.
But the condition $Mf \cap HK = \emptyset$ is equivalent to the condition 
that the endpoint of $p$
does not belong to the $X_1$-component of $v_0$, proving the theorem.

\enddemo

\remark{Remarks}
\roster
\item"1)" The first published proof of the theorem is due to G.A. Niblo \cite{3},
who found a much more simple and elegant argument than the original
 proof of the authors \cite{1}.

\item"2)" A more general result saying that for any finitely generated subgroups

$H_1, \cdots ,H_n$ of a free group $F$  the set 
$H_1 \cdots H_n$ is closed in 
the  profinite topology on $F$ was obtained by L. Ribes and P.A. Zalesskii 
\cite{4},
and by K. Henckell, S.T. Margolis, J.E. Pin and J. Rhodes \cite{5}.

\item"3)"  It is easy to construct a closed subset $B$ of $F$ and a finitely generated
subgroup $K$ such that the product $BK$ is not closed in the profinite 
topology of $F$.
\endroster
\endremark


\Refs
\ref\key 1 
\by R. Gitik and E. Rips
\paper On Separability Properties of Groups
\jour Int J. of Algebra and Computation
\vol 5 \yr 1995 \pages 703-717
\endref  

\ref\key 2 
\by M. Hall, Jr.
\paper Coset Representations in Free Groups
\jour Trans. AMS 
\vol 67 \yr 1949 \pages 431-451
\endref
\ref\key 3 
\by G.A. Niblo 
\paper Separability Properties of Free Groups and Surface 
Groups 
\jour J. of Pure and Applied Algebra 
\vol 78  \yr 1992 \pages 77-84
\endref

\ref\key 4 
\by L. Ribes and P.A. Zalesskii 
\paper On the Profinite Topology on a Free
Group I
\jour  Bull. LMS
\vol 25 \yr 1993  \pages 37-43
\endref

\ref\key 5 
\by K. Henckell, S.T. Margolis, J.E. Pin and J. Rhodes 
\paper Ash's Type II Theorem, Profinite Topology and Malcev Products I
\jour Int J. of Algebra and Computation 
\vol 1  \yr 1991 \pages 411-436
\endref

\ref\key 6 \by J.R. Stallings 
\paper Topology of Finite Graphs 
\jour Invent. Math. 
\vol 71 \yr 1983 \pages 551-565
\endref

\endRefs
\enddocument